%% file: Royset_elastic_rev_arxiv.tex
\documentclass[titlepage,final,10pt]{article}
\input mac.tex

\begin{document}


\begin{center}
\begin{large}
{\bf On Robustness in Nonconvex Optimization with Application to Defense Planning}
\smallskip
\end{large}
\vglue 0.7truecm

{\large Johannes O. Royset}\footnote{Johannes O. Royset, Operations Research Department, Naval Postgraduate School, Monterey, CA93943, USA; joroyset@nps.edu}

\end{center}

\vskip 0.5truecm

\noindent {\bf Abstract}. \quad We estimate the increase in minimum value for a decision that is robust to parameter perturbations as compared to the value of a nominal nonconvex problem. The estimates rely on expressions for subgradients and local Lipschitz moduli of min-value functions and require only the solution of the nominal problem. Across 54 mixed-integer optimization models, the median error in estimating the increase in minimum value is 12\%. The results inform analysts about the possibility of obtaining cost-effective, parameter-robust decisions.

\vskip 0.5truecm

\halign{&\vtop{\parindent=0pt
   \hangindent2.5em\strut#\strut}\cr
{\bf Keywords}: sensitivity analysis, robustness, nonconvex optimization, military operations research.
                         \cr\cr

{\bf Date}:\quad \ \today \cr}

\baselineskip=15pt

\section{Introduction}\label{sec:intro}

An optimization model almost always involves parameters that are unsettled. They may represent unknown probabilities for various demand scenarios, preferences of a risk-averse decision maker, or simply a future operating condition. A practitioner of operations research would need to ensure that decisions recommended by the model is not overly sensitive to changes in such parameters. A senior decision maker would certainly resist the adoption of a decision-support tool that makes recommendations with completely different costs under seemly similar conditions. Likewise, few of us would trust an autonomous system exhibiting an arbitrary ``behavior.'' This brings up the need for sensitivity analysis, a classical subject in optimization dating back at least to the early days of linear programming and the pioneering work by von Neumann on duality theory. Many practical operations research problems lead to nonconvex models due the presence of integer variables, which takes us beyond the classical duality theory of linear and convex programming. Moreover, this theory is traditionally focused on right-hand side perturbations of constraint systems, while many modelers tend to make sure that unsettled parameters appear in an objective function. The reason is widely recognized: constraint systems are often unstable under parametric perturbations \cite{BenkoGfrererOutrata.20,MordukhovichPerezAros.21,Royset.21}. This motivates elastic programming \cite{BrownGraves.75} and models with simple recourse \cite{Wets.74}, where unsettled parameters appear in the objective function.

In this paper, we consider the effect of changing parameters in the objective function of nonconvex optimization problems, but do {\em not} simply quantify the change in minimum values. Instead, we ask the question:
\begin{quote}
  How much more will a decision, robust to parametric changes in the objective function, cost us compared to a decision optimized for nominal parameter values?
\end{quote}
Expressed mathematically in terms of a feasible set $X$ and an objective function $f$ that depends on a parameter vector $u$, the question becomes:
\[
\mbox{Is the minimum value of } ~~~\nnmin_{x\in X} \sup_{u\in U} f(u,x)~~~ \mbox{ much larger than that of } ~~~\nnmin_{x\in X} f(\bar u,x),
\]
where $\bar u$ specifies the nominal parameter values and $U$ is a set of perturbations typically containing $\bar u$. If the answer is no, then there is a cost-effective decision that is robust to parameter changes within the set $U$. If the answer is yes, then robustness comes at a steep price and an optimal decision $x^\star$ for the nominal problem would be severely misleading: there is $\hat u\in U$ such that the cost $f(\hat u,x^\star)$ of that decision is much higher than stipulated by the minimum value of the nominal problem. Thus, the nominal problem is highly sensitive to changes in $u$, which should prompt further study and modeling. For example, one might seek more information about the ``true'' value of $u$.

Using subdifferential calculus \cite{VaAn,Mordukhovich.18,primer}, we estimate the difference between the minimum value of the robust problem and the minimum value of the nominal problem. The resulting estimate is computable from a single solution of the nominal problem. The sensitivity analysis helps us to determine if it is necessary to solve a robust problem, which might entail significant algorithmic development and implementation challenges. We omit a discussion about the multitude of solution approaches to robust optimization problems (often involving duality-based reformulations or decomposition methods) as they arise in semi-infinite programming \cite{Polak.97}, robust convex optimization \cite{BentalElGhaouiNemirovski.09,BertsimasBrownCaramanis.11}, risk-averse decision making \cite{FollmerSchied.11,RockafellarRoyset.15}, distributionally robust optimization \cite{WiesemannKuhnSim.14,RoysetWets.16b}, attacker-defender models \cite{BrownCarlyleSalmeronWood.06,SmithSong.20}, and game theory \cite[Chap. 5, 7, 9]{primer} more generally. This paper centers on how to make a preliminary assessment about whether a cost-efficient, parameter-robust decision is within reach, and omit details about how to eventually solve a robust problem.

The introduction of robustness against parameter variations naturally tends to come with an increase in cost as compared to planning with (assumed) known parameters. This has prompted studies of suitable ``levels of robustness'' as in \cite{BertsimasSim.04,BertsimasGuptaKallus.18}; see \cite{LuShen.21,LiuSaldanhadaGamaWangMao.22} for the contexts of operations management and facility location, respectively. The broad area of machine learning under disturbances to data and/or models \cite{MadryMakelovSchmidtTsiprasVladu.18,SilvaNajafirad.20,WuXiaWang.20,NortonRoyset.22} gives rise to similar trade-offs between robustness and prediction accuracy. Our main theorem on the effect of changes to objective-function parameters in nonconvex optimization is applicable in that setting too. We conjecture that the theorem can help inform analysts about the need for adversarial and/or diametrical training. Still, our present motivation is military operations research and the numerical examples are drawn from that area.

The paper continuous in Section 2 with the main theorem. Section 3 discusses applications from search theory, technology investment for the US Army, and missile load planning for the US Navy. Section 4 provides proofs.

\section{Sensitivity Analysis}

For functions $f_0:\reals^n\to\reals$ and $f_{ik}:\reals^{m_k} \times \reals^n \to \reals$, $i =1, \dots, r_k$, $k = 1, \dots, s$, and a feasible set $X\subset \reals^n$, we consider the {\em nominal problem}
\begin{equation}\label{eqn:modelclass}
\nnmin_{x\in X}  \, f_0(x) + \sum_{k=1}^s \max_{i = 1, \dots, r_k} f_{ik}(\bar u^k,x),
\end{equation}
where the vector $(\bar u^1, \dots, \bar u^s) \in \reals^m$, with $\bar u^k \in \reals^{m_k}$ and $m = \sum_{k=1}^s m_k$, represents nominal values for a parameter vector $u\in\reals^m$ that is unsettled in some sense. The feasible set $X$ is independent of these parameters, but can be nonconvex. The functions $f_0$ and $f_{ik}$ may also be nonconvex, with $f_0$ also possibly nonsmooth.  The problem arises broadly in elastic programming \cite{BrownGraves.75} and goal optimization \cite[Example 4.6]{primer} as well as in the specific defense application models Replenishment At Sea Planner \cite{BrownCarlyle.08,BrownDegrangePriceRowe.17}, Navy Mission Planner \cite{BrownKline.21a}, and Missile Load Planner \cite{BrownKline.21b}. (These models are stated as linear models with additional constraints and slack and surplus variables, but are equivalent to \eqref{eqn:modelclass}.) Section \ref{sec:applic} furnishes additional details and examples.

Given any norm $\|\cdot\|_{(k)}$ on $\reals^{m_k}$, radius $\delta \in \reals$, and center $u\in\reals^m$, we consider the uncertainty set
\[
U(u, \delta) = \big\{ (v^1, \dots, v^s) \in \reals^m ~\big|~ \nmax_{k=1, \dots, s} \|v^k - u^k\|_{(k)} \leq |\delta| \big\}
\]
and the associated {\em robust problem}
\begin{equation}\label{eqn:modelclassRobust}
\nnmin_{x\in X}  \sup_{v\in U(u,\delta)}\, f_0(x) + \sum_{k=1}^s \max_{i = 1, \dots, r_k} f_{ik}(v^k,x).
\end{equation}
If $u = \bar u$ and $\delta = 0$, then the robust problem reduces to the nominal problem. The seemingly unnecessary consideration of negative values of $\delta$ allows us to define the minimum value of the robust problem as a function of $(u,\delta)$ with these arguments being defined on the whole space $\reals^{m+1}$. Formally, the {\em min-value function} corresponding to the robust problem is the function $p:\reals^{m+1} \to [-\infty,\infty]$ defined by
\[
p(u,\delta) = \inf_{x\in X}  \sup_{v\in U(u,\delta)}\, f_0(x) + \sum_{k=1}^s \max_{i = 1, \dots, r_k} f_{ik}(v^k,x).
\]
We seek to quantify how much $p(u,\delta)$ differs from $p(\bar u, 0)$, the minimum value of the nominal problem, for $u$ near $\bar u$ and $\delta$ near zero. Since $p$ is generally nonsmooth and nonconvex, we leverage subdifferential calculus \cite{VaAn,Mordukhovich.18,primer} to estimate this change.

We adopt the following terminology; see for example \cite[Sect. 8.B, 9.A]{VaAn}. For a function $f:\reals^n\to [-\infty, \infty]$ and a point $\bar x$ at which it is finite, we denote the set of subgradients of $f$ at that point by $\partial f(\bar x)$, which then is understood in the general (Mordukhovich, limiting) sense. The function $f$ is locally Lipschitz continuous at $\bar x$ if it is real-valued in a neighborhood of $\bar x$ and there are $\kappa \in [0,\infty)$ and  $\epsilon > 0$ such that
\[
\big|f(x) - f(x')\big| \leq \kappa \|x - x'\|_2 ~~~\mbox{ whenever } \|x - \bar x\|_2 \leq \epsilon, ~\|x' - \bar x\|_2 \leq \epsilon.
\]
It is locally Lipschitz continuous if this holds for all $\bar x\in \reals^n$. The {\em local Lipschitz modulus} of $f$ at $\bar x$ is defined by
\[
\lip f(\bar x) = \limsup_{\substack{x,x'\to \bar x \\ x \neq x'}} \frac{\big|f(x) - f(x')\big|}{\|x-x'\|_2}.
\]
Subgradients and local Lipschitz moduli describe a function locally regardless of its convexity and/or smoothness.

We denote by $\nabla_1 f_{ik}(u^k,x)$ the gradient of $f_{ik}(\cdot,x)$ at the point $u^k$ for a fixed $x$.

\begin{theorem}{\rm (sensitivity).}\label{thm:sen} Suppose that $f_0:\reals^n\to\reals$ is locally Lipschitz continuous, $X\subset \reals^n$ is nonempty and compact, each $f_{ik}:\reals^{m_k} \times \reals^n \to \reals$ is continuously differentiable, $k=1, \dots, s$, $i =1, \dots, r_k$, and $\bar u\in\reals^m$.
Let $M$ be the set of minimizers of \eqref{eqn:modelclass} and, for any $x\in\reals^n$ and $k = 1, \dots, s$, let
\[
I_k(x) = \big\{ i ~\big|~ f_{ik}(\bar u^k,x) = \max_{j=1, \dots, r_k} f_{jk}(\bar u^k, x)\big\}.
\]
Then, the min-value function $p$ is locally Lipschitz continuous at $(\bar u, 0)$ with local Lipschitz modulus
\[
\lip p(\bar u, 0) \leq \max_{\bar x \in M} \sum_{k=1}^s \sum_{i\in I_k(\bar x)} \max_{\|w_i^k\|_{(k)} \leq 1} \sqrt{\big\|\nabla_1 f_{ik}(\bar u^k,\bar x)\big\|_2^2 + \big( \big\langle \nabla_1 f_{ik}(\bar u^k, \bar x), w_i^k \big\rangle\big)^2},
\]
with the right-hand side simplifying to $\sqrt{2} \max_{\bar x \in M} \sum_{k=1}^s \sum_{i\in I_k(\bar x)}\|\nabla_1 f_{ik}(\bar u^k,\bar x)\|_2$ when $\|\cdot\|_{(k)} = \|\cdot\|_2$ for all $k$.

Moreover, for each subgradient $(g,\gamma) \in \partial p(\bar u,0)$, there exist a minimizer $\bar x \in M$ and multipliers $\{w_i^k \in \reals^{m_k}, y^k_i\geq 0, i \in I_k(\bar x), k=1,\dots,s\}$ such that
\[
g = \sum_{k=1}^s \sum_{i \in I_k(\bar x)} y_i^k \nabla_1 f_{ik}(\bar u^k,\bar x), ~~~~~~~~ \gamma = \sum_{k=1}^s \sum_{i \in I_k(\bar x)} y_i^k \big\langle \nabla_1 f_{ik}(\bar u^k,\bar x), w_{i}^k\big\rangle,
\]
where  for all $k$, $\|w_i^k\|_{(k)} \leq 1$, $i=1, \dots, r_k$, and $\sum_{i \in I_k(\bar x)} y_i^k = 1$.

\end{theorem}

The theorem shows that the local Lipschitz modulus of the min-value function $p$ at $(\bar u, 0)$ is bounded by an expression involving the gradient of $f_{ik}$ with respect to its first argument computed at $\bar u^k$ for some worst-case minimizer $\bar x$. While some solvers such as CPLEX \cite{cplex.22} have the capability of identifying all minimizers for mixed-integer linear problems and thus making the full set $M$ available, we expect that most practitioners simple compute the expression for a single $\bar x \in M$ as an estimate of $\lip p(\bar u, 0)$. We see in Section \ref{sec:applic} that such estimates can be quite informative.

The second part of the theorem shows that subgradients of $p$ at $(\bar u,0)$ are of a particular form involving $\nabla_1 f_{ik}(\bar u^k,\bar x)$ and multipliers $y_i^k$ and $w_i^k$. Generally, subgradients provide in some sense more information than the local Lipschitz modulus as they convey directions of more or less growth and not only the magnitude of change.

While the assumption about compactness of $X$ is mild in many operations research models, statistical application may not be naturally expressed using constraints. A close examination of the proof of Theorem \ref{thm:sen} in Section \ref{sec:proofs} shows that the compactness of $X$ can be replaced by an assumption about locally uniform level-boundedness (see \cite[p. 256]{primer}), which in turn is guaranteed when each $f_{ik}$ is bounded from below and $f_0(x)\to\infty$ as $x$ escapes to the horizon. For example, $f_0(x) =\theta\|x\|$, with $\theta>0$ is a possibility.

In the case when the center of the uncertainty set is fixed at $\bar u$ and only the radius $\delta$ changes, we obtain the following specialization of the theorem expressed in terms of the min-value function defined by
\[
q(\delta) = \min_{x\in X}  \max_{v\in U(\bar u,\delta)}\, f_0(x) + \sum_{k=1}^s \max_{i = 1, \dots, r_k} f_{ik}(v^k,x).
\]

\begin{corollary}{\rm (sensitivity to change in radius).}\label{cor:sen}
Suppose that the assumptions in Theorem \ref{thm:sen} hold.
Then, $q$ is locally Lipschitz continuous at $0$ with local Lipschitz modulus
\[
\lip q(0) \leq \max_{\bar x \in M} \sum_{k=1}^s \sum_{i\in I_k(\bar x)} \max_{\|w_i^k\|_{(k)} \leq 1} \big| \big\langle \nabla_1 f_{ik}(\bar u^k, \bar x), w_i^k \big\rangle\big|,
\]
with the right-hand side simplifying to $\max_{\bar x \in M} \sum_{k=1}^s \sum_{i\in I_k(\bar x)}\|\nabla_1 f_{ik}(\bar u^k,\bar x)\|_2$ when $\|\cdot\|_{(k)} = \|\cdot\|_2$ for all $k$.
\end{corollary}

We note that $q$ is univariate and symmetric around the point $0$, which is also a minimizer. Thus, the corollary omits a statement about subgradients of $q$ as they provide no useful information beyond what is already furnished by $\lip q(0)$.

\section{Applications}\label{sec:applic}

To illustrate the theoretical results, we consider three military operations research problems from the areas of search theory, technology investment, and mission planning. The problems result in mixed-integer linear or nonlinear models. Two of the three models are relatively simple to facilitate computations and allow for a complete description. The last model is taken from \cite{BrownKline.21b} and addresses how to load US Navy ships with missiles while accounting for the variety of possible missions the ships might encounter. Throughout, we focus on illustrating Corollary \ref{cor:sen}.

\subsection{Search for Lost Aircraft}

In planning the search for floating debris from an aircraft lost over the ocean, one would divide the area of interest into $s=20$ squares and then determine how much time and resources should be put towards searching each square. Suppose that the search will be carried out by one aircraft flying at $\eta=200$ miles per hour, with a sensor that can reach out $\bar u^k/2$ miles to each side when in square $k$, and that the debris is located in one of the squares and is stationary. Our prior knowledge stipulates the probability that the debris is in square $k$ as $\beta_k>0$, with $\beta_k = 0.02$ for $k = 1$-$5$, $\beta_k = 0.04$ for $k= 6$-$10$, $\beta_k = 0.06$ for $k = 11$-$15$, and $\beta_k = 0.08$ for $k=16$-$20$. The size of each square is $\alpha=3600$ square miles and the total amount of search time available is $\tau=20$ hours. Operational considerations dictate that at most $\kappa$ squares can be searched; we use $\kappa = 8$ or $\kappa = 16$. The problem of determining the time allocation $z = (z_1, \dots, z_{s})$ that minimizes the probability of not finding the debris is addressed by the mixed-integer nonlinear model
\[
\nnmin_{y,z} \,  \sum_{k=1}^s \beta_k \exp( -\eta \bar u^k z_k/\alpha)~~\mbox{ s.t. }~~  \sum_{k=1}^s z_k \leq \tau, ~~\sum_{k=1}^s y_k \leq \kappa, ~~z_k \leq \tau y_k,~~ z_k \geq 0,~~ y_k \in \{0, 1\},~k =1, \dots, s.
\]
We refer to \cite[Sect. 4.2]{StoneRoysetWashburn.16} for further background about search models of this kind.

The model is in the form \eqref{eqn:modelclass} with $x = (y,z)$ and $r_k = 1$ for all $k = 1, \dots, s$. The sensor parameter $\bar u^k$ is unsettled, with a nominal value of $20$ miles. Thus, $f_{1k}(\bar u^k,x) = \beta_k \exp( -\eta \bar u^k z_k/\alpha)$. We consider uncertainty sets with radius $\delta = 5$ miles or $\delta = 10$ miles across all squares (case A), only for squares 7-13 (case B), or only for squares 14-20 (case C). Using Pyomo \cite{Hartetal.17} and the Bonmin solver (with the branch-and-bound option) \cite{bonmin.22}, we compute the minimum value $q(0)$ of the nominal model \eqref{eqn:modelclass} and estimate $\lip q(0)$ in Corollary \ref{cor:sen} using
\[
\widetilde \lip ~q(0) = \sum_{k=1}^s \big |\nabla_1 f_{1k}(\bar u^k,\bar x)\big |,
\]
where $\bar x = (\bar y, \bar z)$ is the obtained minimizer. (We note that these gradients are one-dimensional because $\bar u^k$ is a scalar and they are zero if the parameter is not considered unsettled as occurs for some $k$ in case B and C.) For the two values of $\kappa$ and the two values of $\delta$, Table \ref{tab:search} presents the minimum value of the nominal problem $q(0)$ (column 4), $\widetilde \lip ~q(0)$ (column 5), the estimated minimum value of the robust problem \eqref{eqn:modelclassRobust} as given by
\[
\tilde q(\delta) = q(0) + \delta \,\widetilde \lip ~q(0)
\]
in column 6, the actual minimum value of the robust problem $q(\delta)$ (column 7), and the percentage difference between the estimate $\tilde q(\delta)$ and the actual value $q(\delta)$ (column 8). In this case, the minimum value of the robust problem \eqref{eqn:modelclassRobust} is obtained as easily as that of the nominal problem \eqref{eqn:modelclass} because the worst-case parameter value occurs universally at $\bar u^k - \delta$.

\begin{table}[ht]
\centering
\begin{tabular}{l|r|r|r|r|r|r|r}
case & $\delta$ & $\kappa$ & $q(0)$ & $\widetilde \lip ~q(0)$ & $\tilde q(\delta)$ & $q(\delta)$ & \% error\\
\hline\hline
A   & 5	& 8	&0.448	&0.00389	&0.467	&0.476	&1.8\\
A   & 5	&16	&0.282	&0.01386	&0.351	&0.364	&3.5\\
A   & 10&8	&0.448	&0.00389	&0.486	&0.532	&8.5\\
A   & 10&16	&0.282	&0.01386	&0.421	&0.479	&12.2\\
\hline
B   &5	&8	&0.448	&0.00052	&0.450	&0.451	&0.1\\
B   &5	&16	&0.282	&0.00590	&0.312	&0.317	&1.6\\
B   &10 &8	&0.448	&0.00052	&0.453	&0.457	&0.9\\
B   &10 &16	&0.282	&0.00590	&0.341	&0.365	&6.5\\
\hline
C   &5	&8	&0.448	&0.00336	&0.464	&0.467	&0.5\\
C   &5	&16	&0.282	&0.00688	&0.317	&0.324	&2.3\\
C   &10	&8	&0.448	&0.00336	&0.481	&0.501	&3.8\\
C   &10	&16	&0.282	&0.00688	&0.351	&0.388	&9.4
\end{tabular}
\caption{Minimum values and their estimates in search application.}
\label{tab:search}
\end{table}

We observe that across the cases the estimate $\tilde q(\delta)$ is rather accurate and consistently below the actual minimum value $q(\delta)$. This stems in part from the convexity of the objective function. The accuracy of the estimates deteriorates as the radius is increased from 5 miles to 10 miles, but remains no larger than 12\%.

Returning to the question asked in Section \ref{sec:intro}: is the minimum value of the robust problem much larger than that of the nominal problem? We see from Table \ref{tab:search} that the answer depends on the setting, but  estimates from Corollary \ref{cor:sen} provide a first indication. Robustness comes with little deterioration in performance for case B with $\kappa = 8$ because $q(\delta)$ is only slightly higher than $q(0)$. On the other hand, case A for $\kappa = 16$ results in a significant performance deterioration, which means that an obtained decision from the nominal problem will be rather poor under certain perturbations of the sensor parameter. This might trigger further data collection about the sensor range.

\subsection{Technology Investment}

This example is motivated by a study carried out within the US Army in 2014-2015 to determine the right level of investments in new technologies such as improved nightvision goggles, means to extract water from the air, health sensors, and other futuristic gadgets \cite{TeterRoysetNewman.19}. Specifically, we would like to determine the investments in $T = 10$ technologies to fill needs of the future army pertaining to three areas: efficiency ($a=1$), dominance $(a= 2)$, and expeditionary $(a=3)$.  The goal is to invest as little as possible while satisfying the needs at least to a level of 95\%. Technology $t$ ($t=1, \dots, T$) contributes to area $a$ ($a=1,2,3$) according to the coefficient $u^a_{t}$. Specifically, 
\[
100 - \nsum_{t=1}^{T} u^a_{t} z_t
\]
is the gap in area $a$ under investment $z_t$ in technology $t$, $t = 1, \dots, T$. Thus, we aim for this number to be in the range 0-5, i.e., the gap is between 0 and 5\%.

A difficulty, however, is that the coefficients $u^a_{t}$ are unknown as they relate to future capabilities. To determine these future capabilities, analysts asked subject matter experts from 39 organizations for an assessment; data used in the following is artificial and not representative of the actual survey responses. From the data we produce $N = 100$ possible values (scenarios) for these coefficients: $\bar u_{t}^{a1}, \dots, \bar u_{t}^{aN}$. The problem then becomes to minimize the investment cost plus the expected penalties for failing to have a gap between 0 and 5\%. The penalties are given by the function $h:\reals\to \reals$. Specifically, a gap of $\gamma$ incurs the penalty
\[
h(\gamma) = \max_{i=1,\dots, 5}\alpha_i\gamma + \beta_i,
\]
where $\alpha_1 = -2$ or $-5$, $\beta_1 = 0$, $\alpha_2 = 4$, $\beta_2 = -20$, $\alpha_3 = 40$, $\beta_3 = -380$, $\alpha_4 = 400$, $\beta_4 = -7580$, and $\alpha_5 = \beta_5 = 0$. We see that $h(\gamma)  = 0$ for $\gamma \in [0,5]$, i.e., a zero penalty, but the penalty increases with 4 or more per unit of gap above 5. Negative values of the gap (corresponding to over-investment) are penalized at the rate given by $\alpha_1$.

Investment in technology $t$ should either be zero or between upper and lower bounds $\bar b_t$ and $\underline{b}_t$. This leads to the following model: 
\[
\nnmin_{y\in\reals^{T}, z\in\reals^{T}}  \nsum_{t=1}^{T} z_t  +  \frac{1}{N}\nsum_{n=1}^{N} \nsum_{a=1}^3 h\Big( 100 - \nsum_{t=1}^{T} \bar u_t^{an} z_t \Big) ~\mbox{ s.t. }~ \underline{b}_t y_t \leq z_t \leq \bar b_t y_t, ~y_t \in \{0,1\}, ~t=1, \dots, T.
\]
The problem is of the form \eqref{eqn:modelclass} with $x = (y,z)$, $f_0(x) = \nsum_{t=1}^{T} z_t$, and
\[
f_{ik}(\bar u^k,x) = \alpha_i \Big( 100 - \nsum_{t=1}^{T} \bar u_{t}^{an} z_t \Big) + \beta_i,
\]
where $k = 1, \dots, 3N$ and $\bar u^k$ is the $T$-dimensional vector $(\bar u_{1}^{an}, \dots, \bar u_{T}^{an})$ with $k$ specifying a particular $a$ and $n$. We are concerned about the values of $\bar u^k$ and adopt a distributional robust formulation of the form \eqref{eqn:modelclassRobust} with $\|\cdot\|_{(k)} = \|\cdot\|_2$. We solve the resulting nominal problem as a mixed-integer linear optimization problem using cbc \cite{cbc.22} after introducing auxiliary variables to represent $h$. The robust problem is solved by sampling 500 points within the uncertainty set, expanding the mixed-integer linear model accordingly, and applying cbc. (The sampling approach furnishes a lower bound on the minimum value of the robust problem, which we compare with an accurate objective function evaluation. When reporting estimation accuracy below, we report the least favorable number produced by this evaluation and the lower bound.)

We consider radius of size $\delta = 0.05, 0.10$, and $0.20$, which correspond to roughly 5\%, 10\%, and 20\% change away from $\bar u_t^{an}$; their values tend to be around 1. We also consider two values of $\alpha_1$ and three versions of the bounds $\underline{b}_t$ and $\bar b_t$: base, low upper bound, and high lower bound. The results are summarized in Table \ref{tab:army} using the same notation as in Table \ref{tab:search}. The second to last column reports an interval of lower and upper bound on the minimum value of the robust problem due our inexact solution of that problem.

\begin{table}[ht]
\centering
\begin{tabular}{l|r|r|r|c|r|c|r}
bound & $\alpha_1$ & $\delta$ & $q(0)$ & $\widetilde \lip ~q(0)$ & $\tilde q(\delta)$ & $q(\delta)$ & \% error\\
\hline\hline
base      & $-2$ & 0.05	&130	&379	&149	&   [156.5  156.9]	&	4.5\\
base      & $-2$ & 0.1	&130	&379	&168	&	[185.1  186.3]	&   9.4\\
base      & $-2$ & 0.2	&130	&379	&206	&	[249.8	252.0]&	18.1\\
\hline
low upper & $-2$ & 0.05	&138	&427	&159	&	[164.4	164.8]&	3.0\\
low upper & $-2$ & 0.1	&138	&427	&181	&	[191.7	193.0]&	6.1\\
low upper & $-2$ & 0.2	&138	&427	&223	&	[250.8	253.1]&	11.6\\
\hline
high lower& $-2$ &0.05	&151	&204	&161	&	[168.2	168.6]	&   4.0\\
high lower& $-2$ &0.1	&151    &204	&172	&	[195.4	196.3] &	12.4\\
high lower& $-2$ &0.2	&151	&204	&192	&	[265.9	268.7&	28.6\\
\hline
base      & $-5$ & 0.05	&141	&562	&169	&	[175.8	176.5]&	4.2\\
base      & $-5$ & 0.1	&141	&562	&197	&	[224.3	226.2]&	12.9\\
base      & $-5$ & 0.2	&141	&562	&253	&	[340.9	345.7]&	27.0\\
\hline
low upper & $-5$ & 0.05	&151	&533	&177	&	[184.3	184.8]&	3.8\\
low upper & $-5$ & 0.1	&151	&533	&204	&	[226.9	228.5]&	10.5\\
low upper & $-5$ & 0.2	&151	&533	&257	&	[340.4	345.8]&	25.8\\
\hline
high lower& $-5$ &0.05	&154	&300	&169	&	[183.2	184.0]&	8.0\\
high lower& $-5$ &0.1	&154	&300	&184	&	[236.2	238.4]&	22.9\\
high lower& $-5$ &0.2	&154	&300	&214	&	[368.5	374.0]&	43.3\\
\end{tabular}
\caption{Minimum values and their estimates in investment application.}
\label{tab:army}
\end{table}

Table \ref{tab:army} shows that the estimate $\tilde q(\delta)$ is moderately accurate and again consistently below the actual minimum value $q(\delta)$. The accuracy of the estimates worsens for larger $\delta$ as one could expect. Moreover, the case with high lower bounds on the investments seems easier to robustify; the increase from $q(0)$ to $q(\delta)$ is smaller than in the other cases. Thus, we have identified a case where a distributionally robust decision can be achieved at a relatively low cost. The base case appears more challenging as a robust decision tends to increase the cost significantly. We can reach these assessments by simply solving the nominal problem and leveraging Corollary \ref{cor:sen}.

\subsection{Missile Load Planner}

This example is taken from \cite{BrownKline.21b}, which develops a  mixed-integer optimization model referred to as the Missile Load Planner. Given a set of possible missions, warplans, and deployment cycles, the model assigns missiles of many different types from a finite inventory to a group of ships so that several concerns are addressed: The number of missiles assigned to a ship is within its storage capacity; the mix of missiles matches or exceeds the requirements of the selected missions, which are also assigned different priorities; the number of ships is sufficient for each selected mission; and the change from an existing load onboard a ship is not excessive. In summary, the model takes the following form:\\

{\bf Indices}

$w$ ~warplan

$m$ ~mission

$d$ ~~deployment cycles

$s$ ~~\,ship

$y,y'$ missile type\\

{\bf Data}

$priority_m$: ~priority of mission $m$

$shipsreq_m$: ~ships required by mission $m$

$shipsshortpen_m$: ~ship shortfall penalty for mission $m$

$missilesmin_{m,y}$: ~minimum number of missiles of type $y$ on each ship for mission $m$

$missileshortpen_{m,y}$: ~missile shortfall penalty for mission $m$, type $y$

$altmissilepen_{m,y,y'}$: ~penalty for substituting missile type $y'$ for $y$ in mission $m$

$changepen$:~ penalty for adjusting prior loadout\\

{\bf Decision Variables}

$ASSIGN_{w,d,m,s}$: ~assign ship $s$ to warplan $w$, deployment cycle $d$, mission $m$ [binary]

$MISSION_{w,d,m}$: ~commit to warplan $w$, cycle $d$, mission $m$ [binary]

$COMMIT_{w,d,m,s,y,y'}$: ~warplan $w$, deployment cycle $d$, mission $m$, ship $s$, want type $y$, commit type $y'$

$CHANGE_{s,y}$: ~number of $y$ missiles changed in ship $s$\\

{\bf Formulation}

\begin{align*}
&\nnmin -\sum_{w,d,m} priority_m MISSION_{w,d,m} + \sum_{w,d,m,s,y,y'} altmissilepen_{m,y,y'} COMMIT_{w,d,m,s,y,y'}\\
&+ \sum_{w,d,m} \max\bigg\{0, ~shipsshortpen_m \bigg(shipsreq_m MISSION_{w,d,m} - \sum_s ASSIGN_{w,d,m,s} \bigg)\bigg\}\\
&+ \sum_{s,y} changepen \,CHANGE_{s,y}\\
&+ \sum_{w,d,m,y} \max\bigg\{0, ~missilesshortpen_{m,y} \bigg(missilesmin_{m,y} MISSION_{w,d,m} - \sum_{s,y'} COMMIT_{w,d,m,s,y,y'}  \bigg)  \bigg\}
\end{align*}

subject to constraints listed as (2), (3), (5)-(7), (9), (10), (12)-(17) in \cite{BrownKline.21b}.\\

The first term in the objective function assigns ``rewards'' for executing missions. The second term specifies a penalty for substituting a missile type for a less ideal one. The third term assigns a penalty if too few ships are assigned to a mission. The fourth term invokes a penalty when a legacy loadout of missiles is changed. The last term assigns a penalty when fewer than $missilesmin_{m,y}$ of type $y$ is committed to mission $m$. We refer to \cite{BrownKline.21b} for further details and explanations. After the introduction of auxiliary variables to linearize the max-expressions, a typical model instance has 3,435 constraints and 6,756 variables of which 3,298 are integer variables.

The minimum number of missiles of type $y$ in mission $m$, $missilesmin_{m,y}$, is unsettled and plays the role of $\bar u^k$ in \eqref{eqn:modelclass}. The first four terms in the objective function give rise to $f_0$ in \eqref{eqn:modelclass}, while the last term furnishes the sum in \eqref{eqn:modelclass}. Thus, in the present example, the max-expression in \eqref{eqn:modelclass} is always over two elements ($r_k = 2$), with $f_{1k}$ simply being the zero function and $f_{2k}$ being given by
\[
missilesshortpen_{m,y} \bigg(missilesmin_{m,y} MISSION_{w,d,m} - \sum_{s,y'} COMMIT_{w,d,m,s,y,y'}  \bigg),
\]
which indeed depends on the unsettled parameter $missilesmin_{m,y}$. We consider an uncertainty set $U(\bar u, \delta)$ with $\delta = 1$. Since $\bar u^k$ is a scalar, this means that the perturbation entails changing any combination of the parameters $missilesmin_{m,y}$ by one missile. (While the nominal values of these parameters can be as high as 100, it can also be as low as 0. A change of one missile is thus sizable.)

We solve the nominal problem \eqref{eqn:modelclass}  to optimality using CPLEX \cite{cplex.22}. In this case, the robust problem \eqref{eqn:modelclassRobust} is solved in the same manner as the characteristics of the model allow us to deduce apriori that the worst case involves {\em higher} values of $missilesmin_{m,y}$.

Table \ref{tab:load} shows results parallel to those of Table \ref{tab:search} across 24 cases. We consider three cases of $missilesmin_{m,y}$, i.e., $\bar u^k$ in \eqref{eqn:modelclass}, as indicated in column 1 of the table. These corresponds to low, medium, and high missile requirements as encoded in $missilesmin_{m,y}$. The penalty for changing a legacy loadout is varied according to column 2. Column 3 specifies whether ships used in a deployment cycle of a warplan must be used in all deployment cycles of that plan; see constraints (14)-(17) in \cite{BrownKline.21b}.

\begin{table}[ht]
\centering
\begin{tabular}{c|c|c|r|c|r|c|r}
$missilesmin_{m,y}$ & $changepen$ & cycle constr. & $q(0)$ & $\widetilde \lip ~q(0)$ & $\tilde q(\delta)$ & $q(\delta)$ & \% error\\
\hline\hline
low & 0	&inactive	&1584	&5350	&6934	&1968	&313\\
low &0	&active	&1584	&5180	&6764	&1968	&302\\
low &5	&inactive	&3772	&3780	&7552	&4238	&87\\
low &5	&active	&4138	&3980	&8118	&4579	&85\\
low &10	&inactive	&5347	&3000	&8347	&5809	&47\\
low &10	&active	&5970	&3630	&9600	&6464	&52\\
low &20	&inactive	&7908	&2450	&10358	&8405	&24\\
low &20	&active	&8959	&1850	&10809	&9512	&14\\
\hline
medium &0	&inactive	&4414	&4150	&8564	&4860	&83\\
medium &0	&active	&4414*	&4800	&9214	&4899	&97\\
medium &5	&inactive	&6950	&2770	&9720	&7358	&33\\
medium &5	&active	&7400	&2620	&10020	&7862	&29\\
medium &10	&inactive	&8530	&2580	&11110	&9006	&24\\
medium &10	&active	&9469	&1010	&10479	&9899	&6\\
medium &20	&inactive	&11463	&2700	&14163	&12021	&18\\
medium &20	&active	&12562	&1160	&13722	&13094	&4\\
\hline
high &0	&inactive	&6853	&4320	&11173	&7487	&53\\
high &0	&active	&6897**	&4040	&10937	&7646	&47\\
high &5	&inactive	&9665	&2100	&11765	&10158	&16\\
high &5	&active	&10341	&1860	&12201	&10862	&12\\
high &10	&inactive	&11535	&2350	&13885	&12069	&15\\
high &10	&active	&12438	&1160	&13598	&13045	&4\\
high &20	&inactive	&15025	&2580	&17605	&15650	&13\\
high &20	&active	&16058	&1080	&17138	&16650	&3
\end{tabular}
\caption{Minimum values and their estimates in Missile Load Planner. (*Solved to 0.68\% tolerance after one hour. **Solved to 0.11\% tolerance after one hour.)}
\label{tab:load}
\end{table}

Since the problem is highly combinatorial, the effect of parameter changes becomes less predictable. The minimum-value estimates from Corollary \ref{cor:sen} (see column 6) exceed the true values, sometimes significantly. Still, the prediction error is often within 20\% for cases with higher minimum missile requirements; see the last eight rows of Table \ref{tab:load}. We can use the estimates to assess whether a robust loadout plan with a low objective function value is achievable. For example, in the last row of Table \ref{tab:load} there is a relatively small change from $q(0)$ to $\tilde q(\delta)$ indicating the possibility of obtaining a ``good'' loadout plan that accounts for ambiguity about $missilesmin_{m,y}$. This is indeed possible because $q(\delta)$ is only slightly larger than $q(0)$. Switching to the first case with a medium missile requirement, we see that Corollary \ref{cor:sen} predicts a relatively large increase from $q(0) = 4414$ to $\tilde q(\delta) = 8564$. While the actually increase turns out to be more moderate $(q(\delta) = 4860)$, the percent increase is three times as large as that of the last-row case. Thus, the estimate from Corollary \ref{cor:sen} rightly flags the first case with medium missile requirement as potentially concerning: robustness to changes in $missilesmin_{m,y}$ comes with a significant worsening of the performance as measured by the objective function value.

\section{Proofs}\label{sec:proofs}

\state Proof of Theorem \ref{thm:sen}. For any $x\in\reals^n$, $u=(u^1, \dots, u^s)\in\reals^m$, and $\delta\in\reals$, the maximization over $U(u,\delta)$ can be moved ``inside'' the sum, i.e.,
\[
\max_{v\in U(u,\delta)}\, f_0(x) + \sum_{k=1}^s \max_{i = 1, \dots, r_k} f_{ik}(v^k,x) = f_0(x) + \sum_{k=1}^s \max_{i = 1, \dots, r_k} \max_{\|w^k_i\|_{(k)} \leq 1} f_{ik}(u^k + \delta w_i^k, x),
\]
with the maxima being attained because $f_{ik}(\cdot,x)$ is continuous and $U(u,\delta)$ is nonempty and compact.
Let $\iota_X(x) = 0$ if $x\in X$ and $\iota_X(x) = \infty$ otherwise. Then, $p(u,\delta) = \min_{x\in \reals^n} \phi(u,\delta,x)$, where
\begin{align*}
\phi(u,\delta,x) & = \iota_X(x) + f_0(x) + \sum_{k=1}^s \psi_k(u^k,\delta,x)\\
\psi_k(u^k,\delta,x) & = \max_{i = 1, \dots, r_k} \psi_{ik}(u^k,\delta,x), ~~k=1, \dots, s\\
\psi_{ik}(u^k,\delta,x) & = \max_{\|w_i^k\|_{(k)} \leq 1} f_{ik}(u^k + \delta w_i^k, x), ~~k=1, \dots, s, ~i = 1, \dots, r_k.
\end{align*}
The function $\psi_{ik}$ is continuous by \cite[Proposition 6.29]{primer} and also real-valued; see \cite[Theorem 4.9]{primer}. The function $\psi_k$ is continuous and real-valued, which we deduce from the proof of \cite[Proposition 4.66]{primer}. The function $f_0$ is continuous and real-valued by assumption. Since $X$ is closed and nonempty, we conclude that $\phi$ is lower semicontinuous and proper (see p.187 and p.73 in \cite{primer}, respectively, for definitions). Thus, we can bring in \cite[Theorem 10.13]{VaAn}; the requirement of level-boundedness in that theorem holds because $X$ is compact. This theorem directs us to compute $\partial \psi(\bar u, 0,\bar x)$ and the set of horizon subgradients $\partial^\infty \psi(\bar u, 0,\bar x)$ at a minimizer $\bar x$ of the nominal problem \eqref{eqn:modelclass}; see \cite[Definition 4.60]{primer}. Actually, computing outer approximations suffice.

By \cite[Proposition 6.30]{primer}, $\psi_{ik}$ is locally Lipschitz continuous. The same holds for $\psi_k$ by \cite[Proposition 4.68]{primer}. For $u\in\reals^m$, $\delta\in \reals$, and $x\in \reals^n$, we utilize \cite[Example 4.70]{primer} to obtain the outer approximation
\begin{equation}\label{eqn:partialeq}
\partial \phi(u,\delta,x) \subset \{0\}^{m+1} \times N_X(x) + \{0\}^{m+1} \times \partial f_0(x) + \sum_{k=1}^s \partial \psi_k(u^k,\delta,x),
\end{equation}
where $N_X(x) = \partial \iota_X(x)$ is the normal cone of $X$ at $x$; see for example \cite[Section 4.G]{primer}. Next, we consider $\partial \psi_k(u^k,\delta,x)$ appearing at the end of \eqref{eqn:partialeq}. The function $\psi_k$ can be written as a composition:  $\psi_k(u^k,\delta,x) = h_k(F_k(u^k,\delta,x))$, where $h_k(z) = \nmax_{i=1, \dots, r_k} z_i$ and $F_k:\reals^{m_k}\times \reals\times\reals^n \to \reals^{r_k}$ with
\[
F_k(u^k,\delta, x) = \big( \psi_{1k}(u^k,\delta,x), \dots, \psi_{r_kk}(u^k,\delta,x)  \big).
\]
The mapping $F_k$ is locally Lipschitz continuous because each $\psi_{ik}$ has that property. For $y\in \reals^{m_k}$, let $f_k(u^k,\delta,x,y)$ $=$ $\langle F_k(u^k,\delta,x), y\rangle$. Since the set of horizon subgradients $\partial^\infty h_k(z) = \{0\}$ because $h_k$ is real-valued and convex (see \cite[Proposition 4.65]{primer}), the qualification (6.15) in \cite[Theorem 6.23]{primer} holds. Thus, by that theorem,
\begin{equation}\label{eqn:psik}
\partial \psi_k(u^k,\delta, x) = \partial (h_k \circ F_k)(u^k,\delta, x) \subset \bigcup_{y\in \partial h_k (F_k(u^k,\delta,x))}  \partial_{-y} f_k(u^k,\delta,x,y),
\end{equation}
where $\partial_{-y}$ indicates that subgradients are computed with respect to all variables expect $y$.
Now $y\in \partial h_k(z)$ if and only if
\[
\sum_{i=1}^{r_k} y_i = 1, ~~~~~~ y_i \geq 0 ~\mbox{ if } ~z_i = h_k(z) ~~\mbox{ and otherwise }~ y_i = 0, ~~i=1, \dots, r_k
\]
as can be seen from \cite[Proposition 4.66]{primer} and its proof. Thus, it suffices to consider nonnegative $y_i$ in the expression
\[
f_k(u^k,\delta,x,y) = \big\langle F_k(u^k,\delta,x), y\big\rangle = \sum_{i=1}^{r_k}  y_i \psi_{ik}(u^k,\delta,x).
\]
We note that the function $(u^k,\delta,x)\mapsto y_i \psi_{ik}(u^k,\delta,x)$ is locally Lipschitz continuous because $\psi_{ik}$ has that property. Thus, \cite[Example 4.70]{primer} applies and yields
\begin{equation}\label{eqn:fk}
\partial_{-y} f_k(u^k,\delta,x,y) \subset \sum_{i=1}^{r_k}  \partial (y_i \psi_{ik})(u^k,\delta,x).
\end{equation}
If $y_i = 0$, then $\partial (y_i \psi_{ik})(u^k,\delta,x) = \{0\}$. If $y_i > 0$, then by \cite[Eq. (4.13)]{primer} $\partial (y_i \psi_{ik})(u^k,\delta,x) = y_i \partial \psi_{ik}(u^k,\delta,x)$.

Next, we utilize \cite[Proposition 6.30]{primer} to compute
\[
\partial \psi_{ik}(u^k,\delta,x) = \con\Big\{ \nabla \tilde f_{ik}(u^k,\delta,x;\bar w_i^k) ~\Big|~ \bar w_i^k \in \nargmax_{w}\big\{ f_{ik}(u^k + \delta w,x) \big| \|w\|_{(k)} \leq 1 \big\}\Big\},
\]
where $\con S$ is the convex hull of a set $S$ and, for $w\in \reals^{m_k}$, we adopt the notation $\tilde f_{ik}(u^k,\delta,x; w ) =  f_{ik}(u^k + \delta w, x)$. Thus, $\nabla \tilde f_{ik}(u^k,\delta,x; w)$ is the gradient of $\tilde f_{ik}(\cdot,\cdot,\cdot; w )$ at $(u^k,\delta,x)$, i.e.,
\[
\nabla \tilde f_{ik}(u^k,\delta,x;w) = \begin{bmatrix}
\nabla_1 f_{ik}(u^k + \delta w, x)\\
   \big\langle \nabla_1 f_{ik}(u^k + \delta w, x), w\big\rangle\\
   \nabla_2 f_{ik}(u^k + \delta w, x)
\end{bmatrix}.
\]
Here, $\nabla_1$ and $\nabla_2$ indicate differentiation with respect to the first and second argument, respectively. Then,
\[
\partial (y_i \psi_{ik})(u^k,\delta,x) = y_i \con \left.\left\{ \begin{bmatrix}
\nabla_1 f_{ik}(u^k + \delta \bar w_i^k, x)\\
   \big\langle \nabla_1 f_{ik}(u^k + \delta \bar w_i^k, x), \bar w_i^k\big\rangle\\
   \nabla_2 f_{ik}(u^k + \delta \bar w_i^k, x)
\end{bmatrix} \right|~ \bar w_i^k \in \nargmax_{w}\big\{ f_{ik}(u^k + \delta w,x) \big| \|w\|_{(k)} \leq 1 \big\} \right\}.
\]
This expression can be inserted into \eqref{eqn:partialeq} via \eqref{eqn:psik} and \eqref{eqn:fk}. We obtain the outer approximation
\[
  \partial \phi(u,\delta,x) \subset \{0\}^{m+1} \times N_X(x) + \{0\}^{m+1} \times \partial f_0(x)\\
  + \sum_{k=1}^s \bigcup_{y^k \in Y_k(u^k,\delta,x)} \sum_{i = 1}^{r_k} \partial (y_i^k \psi_{ik})(u^k,\delta,x),
\]
where
\[
Y_k(u^k,\delta,x) = \Big\{ y \in [0,\infty)^{r_k}  ~\Big|~ \nsum_{i=1}^{r_k} y_i = 1, ~y_i = 0 \mbox{ if } \psi_{ik}(u^k,\delta,x) < \psi_k(u^k,\delta,x)   \Big\}.
\]
Evaluation these subgradients at the point $(\bar u,0,x)$, we obtain the simplification
\begin{align*}
  \partial \phi(\bar u,0,x) & \subset \{0\}^{m+1} \times N_X(x) + \{0\}^{m+1} \times \partial f_0(x)\\
  &  ~~~~~~~ + \sum_{k=1}^s \bigcup_{y^k \in Y_k(\bar u^k,0,x)} \sum_{i = 1}^{r_k} y_i^k \con \left.\left\{ \begin{bmatrix}
\nabla_1 f_{ik}(\bar u^k, x)\\
   \big\langle \nabla_1 f_{ik}(\bar u^k, x), \bar w^k_i\big\rangle\\
   \nabla_2 f_{ik}(\bar u^k, x)
\end{bmatrix} \right|~ \|\bar w_i^k\|_{(k)} \leq 1 \right\}.
\end{align*}

Suppose that $(g,\gamma) \in \partial p(\bar u,0)$. By \cite[Theorem 10.13]{VaAn}, there is a minimizer $\bar x$ of \eqref{eqn:modelclass} such that $(g,\gamma,0) \in \ \partial \phi(\bar u,0,\bar x)$. Using the above inclusion and Caratheodory's theorem (see for example \cite[p. 203]{primer}), we obtain that there are $\{y^k\in Y_k(\bar u^k,0,\bar x), k=1, \dots, s\}$ and $\{w_{ij}^k \in \reals^{m_k}, \lambda_{ij}^k\geq 0, i = 1, \dots, r_k, j=0, 1, \dots, m_k + 1 + n, k=1, \dots, s\}$ with $\sum_{j=0}^{m_k + 1+n} \lambda_{ij}^k = 1$ and $\|w_{ij}^k\|_{(k)} \leq 1$ for all $ i = 1, \dots, r_k$, $j=0, 1, \dots, m_k + 1 + n$, and $k = 1,\dots, s$, such that
\begin{align*}
g & = \sum_{k=1}^s \sum_{i=1}^{r_k} y_i^k \sum_{j=0}^{m_k + 1+n} \lambda_{ij}^k \nabla_1 f_{ik}(\bar u^k,\bar x)\\
\gamma & = \sum_{k=1}^s \sum_{i=1}^{r_k} y_i^k \sum_{j=0}^{m_k + 1+n} \lambda_{ij}^k  \big\langle \nabla_1 f_{ik}(\bar u^k,\bar x), w_{ij}^k\big\rangle\\
0 & \in N_X(\bar x) + \partial f_0(\bar x) + \sum_{k=1}^s \sum_{i=1}^{r_k} y_i^k \sum_{j=0}^{m_k + 1+n} \lambda_{ij}^k \nabla_2 f_{ik}(\bar u^k,\bar x).
\end{align*}
Let $w_i^k = \sum_{j=0}^{m_k + 1+n} \lambda_{ij}^k w_{ij}^k$. Since $\sum_{j=0}^{m_k + 1+n} \lambda_{ij}^k = 1$, we find that  $\|w_i^k\|_{(k)} \leq 1$ and
\begin{align*}
g & = \sum_{k=1}^s \sum_{i=1}^{r_k} y_i^k \nabla_1 f_{ik}(\bar u^k,\bar x)\\
\gamma & = \sum_{k=1}^s \sum_{i=1}^{r_k} y_i^k \big\langle \nabla_1 f_{ik}(\bar u^k,\bar x), w_{i}^k\big\rangle\\
0 & \in N_X(\bar x) + \partial f_0(\bar x) + \sum_{k=1}^s \sum_{i=1}^{r_k} y_i^k \nabla_2 f_{ik}(\bar u^k,\bar x).
\end{align*}
This confirms the assertion about subgradients.

To establish the assertion about local Lipschitz continuity of $p$ at $(\bar u, 0)$ and the associated local Lipschitz modulus, we invoke \cite[Corollary 10.14]{VaAn}, which applies when, regardless of minimizer $\bar x$ of \eqref{eqn:modelclass}, one has
\[
(\bar g,\bar\gamma,0) \in \partial^\infty \phi(\bar u,0,\bar x) ~~~\Longrightarrow~~~ \bar g =0, ~~ \bar \gamma = 0.
\]
From \cite[Example 4.70]{primer}, we confirm that $\partial^\infty \phi(\bar u,0,\bar x) \subset \{0\}^{m+1} \times N_X(\bar x)$. Thus, the corollary applies and the conclusion follows.\eop\\

\state Proof of Corollary \ref{cor:sen}. One can argue as in the proof of Theorem \ref{thm:sen} but account for the simplifications stemming from the fact that $u$ is fixed at $\bar u$.\eop\\

\state Acknowledgement. The author is grateful to Gerald Brown (NPS) for providing an implementation of the Missile Load Planner in GAMS. This work is supported in part by the Office of Naval Research (Mathematical and Resource Optimization Program) under MIPR N0001422WX01445.

\bibliographystyle{plain}
\bibliography{refs}

\end{document}

%% file: mac.tex
\usepackage{theorem}
\usepackage{enumerate}
\usepackage{array}
\usepackage[reqno]{amsmath}
\usepackage{amssymb}
\usepackage{mathtools}
\usepackage{latexsym}
\usepackage{makeidx}
\usepackage{fancybox}
\input epsf.sty
\usepackage[dvips]{graphicx,color}
\usepackage{showlabels}
\usepackage{subcaption}

\theoremstyle{change}
\newtheorem{proclaim}{PROCLAIM}[section]
\newtheorem{theorem}[proclaim]{Theorem}

\newtheorem{corollary}[proclaim]{Corollary}

\numberwithin{equation}{section}

\outer\def\proclaim #1. #2\par{\medbreak \noindent{\bf#1.\enspace}{\sl#2}\par
  \ifdim\lastskip<\medskipamount
  \removelastskip\penalty55\medskip\fi}
\def\state #1. { \noindent{\bf#1.\enspace}}
\def\algo #1. { \noindent{\bf#1.\enspace}}

\DeclareMathOperator{\con}{con}

\DeclareMathOperator{\lip}{lip}

\newcommand{\comp}{\,{\raise 1pt \hbox{$\scriptstyle\circ$}}\,}

\newcommand{\reals}{\mathbb{R}}

\newcommand{\natnums}{{{\rm l} \kern -.13em {\rm N} }}

\newcommand{\snats}{{I\kern -.29em N}}
\newcommand{\rats}{{Q\kern -.64em \raise 1pt \hbox{$\scriptstyle |$}\;\,}}
\newcommand{\srats}
	{{Q\kern -.56em \raise 1.2pt \hbox{$\scriptscriptstyle /$}\,}}
\newcommand{\ints}{Z\kern -.46em Z}

\newcommand{\pluss}{\hskip1pt \raise1pt\vbox{\hrule width6pt \vskip1pt \hrule
                    width6pt} \kern-4pt{\lower1pt\hbox{\vrule height6pt
		    \kern1pt\vrule height6pt}}\hskip5pt}
\newcommand{\eop}
	{\hfill{$\vcenter{\hrule height1pt \hbox{\vrule width1pt height5pt
   	 \kern5pt \vrule width1pt} \hrule height1pt}$} \medskip}

\newcommand{\setd}{{ d \kern -.15em l}}
\newcommand{\hatsetd}{ d \hat{\kern -.15em l }}

\renewcommand{\epsilon}{\varepsilon}
\renewcommand{\phi}{\varphi}

\newcommand{\tto}{\;{\lower 1pt \hbox{$\rightarrow$}}\kern -12pt
           \hbox{\raise 2.5pt \hbox{$\rightarrow$}}\;}
\newcommand{\overto}[1]{\,{\raise 0pt\hbox{$\rightarrow$}}\kern -9pt
     \hbox{\lower 3pt \hbox{$\scriptscriptstyle#1$}}\hskip6pt}
\newcommand{\underto}[1]{\,{\lower 1pt\hbox{$\rightarrow$}}\kern -9pt
     \hbox{\raise 4pt \hbox{$\,\scriptscriptstyle#1$}}\hskip7pt}
\newcommand{\bigoverto}[1]{{\raise 0pt\hbox{$\,\longrightarrow$}}\kern -16pt
     \hbox{\lower 3pt \hbox{$\scriptscriptstyle#1$}}\hskip4pt}
\newcommand{\bigunderto}[1]{\,{\lower 1pt\hbox{$\longrightarrow$}}\kern -16pt
     \hbox{\raise 4pt \hbox{$\,\scriptscriptstyle#1$}}\hskip6pt}
\newcommand{\bigbigto}[2]{\,{\raise 0pt\hbox{$\,\longrightarrow$}}\kern -16pt
     \hbox{\lower 3pt \hbox{$\scriptscriptstyle#2$}}\kern -10pt
     \hbox{\raise 4pt \hbox{$\,\scriptscriptstyle#1$}}\hskip7pt}
\newcommand{\downto}{{\raise 1pt \hbox{$\scriptscriptstyle \,\searrow\,$}}}
\newcommand{\upto}{{\raise 1pt \hbox{$\scriptscriptstyle \,\nearrow\,$}}}

\newcommand{\notimply}
	{\quad\hbox{$\Longrightarrow \kern -14pt {/}$}\hskip6pt\quad}

\newcommand{\lto}{\,{\lower 1pt\hbox{$\rightarrow$}}\kern -10pt
     \hbox{\raise 4pt \hbox{$\, \scriptstyle l$}}\hskip7pt}
\newcommand{\eto}{\,{\lower 1pt\hbox{$\rightarrow$}}\kern -11pt
     \hbox{\raise 4pt \hbox{$\, \scriptstyle e$}}\hskip7pt}
\newcommand{\hto}{\,{\lower 1pt\hbox{$\rightarrow$}}\kern -11pt
     \hbox{\raise 4pt \hbox{$\, \scriptstyle h$}}\hskip7pt}
\newcommand{\pto}{\,{\lower 1pt\hbox{$\rightarrow$}}\kern -11pt
     \hbox{\raise 4.5pt \hbox{$\, \scriptstyle p$}}\hskip7pt}
\newcommand{\cto}{\,{\lower 1pt\hbox{$\rightarrow$}}\kern -11pt
     \hbox{\raise 4pt \hbox{$\, \scriptstyle c$}}\hskip7pt}
\newcommand{\gto}{\,{\lower 1pt\hbox{$\rightarrow$}}\kern -11pt
     \hbox{\raise 4.5pt \hbox{$\, \scriptstyle g$}}\hskip7pt}
\newcommand{\sto}{\,{\lower 1pt\hbox{$\rightarrow$}}\kern -11pt
     \hbox{\raise 4pt \hbox{$\, \scriptstyle s$}}\hskip7pt}
\newcommand{\awto}{\,{\lower 1pt\hbox{$\rightarrow$}}\kern -15pt
     \hbox{\raise 4pt \hbox{$\, \scriptstyle aw$}}\hskip7pt}
\def\Nto{\,{\raise 1pt\hbox{$\rightarrow$}}\kern -13pt
     \hbox{\lower 3pt \hbox{$\, \scriptstyle N$}}\hskip7pt}
\def\Cto{\,{\raise 1pt\hbox{$\rightarrow$}}\kern -14pt
     \hbox{\lower 3pt \hbox{$\, \scriptstyle C$}}\hskip7pt}
\def\fto{\,{\raise 1pt\hbox{$\rightarrow$}}\kern -14pt
     \hbox{\lower 3pt \hbox{$\, \scriptstyle f$}}\hskip7pt}

\newcommand{\low}[1]{{\lower1pt \hbox{$\scriptstyle #1$}}}
\newcommand{\loww}[1]{{\lower2pt \hbox{$\scriptstyle #1$}}}
\newcommand{\high}[1]{{\raise1pt \hbox{$\scriptstyle #1$}}}

\newcommand{\nsum}{\mathop{\sum}\nolimits}

\newcommand{\nmax}{\mathop{\rm max}\nolimits}

\newcommand{\nnmin}{\mathop{\rm minimize}}

\newcommand{\nargmax}{\mathop{\rm argmax}\nolimits}

\newcommand{\lwdy}[2]{\mathrel{\mathop
        {\raisebox{0.1ex}{\null$#1$}}{\hbox{\kern -1.0em
	{\raisebox{-0.8ex}{$\scriptstyle{\;\to #2}$}}}}}}
\newcommand{\lwwdy}[2]{\mathrel{\mathop
        {\raisebox{0.2ex}{\null$#1$}}{\hbox{\kern -1.0em
	{\raisebox{-1.1ex}{$\scriptstyle{\;\to #2}$}}}}}}
\newcommand{\slwdy}[2]{\scriptsize{{\mathrel{\mathop
        {\raisebox{0.1ex}{\null$#1$}}{\hbox{\kern -1.0em
	{\raisebox{-0.8ex}{$\scriptstyle{\;\to #2}$}}}}}}}}
\newcommand{\slwwdy}[2]{\scriptsize{{\mathrel{\mathop
        {\raisebox{0.2ex}{\null$#1$}}{\hbox{\kern -1.0em
	{\raisebox{-1.1ex}{$\scriptstyle{\;\to #2}$}}}}}}}}

\definecolor{lightgray}{gray}{0.75}
\definecolor{myred}{rgb}{0.55,0,0}
\definecolor{myblue}{rgb}{0,0,0.5} 
\definecolor{mygreen}{rgb}{0,0.5,0} 
\definecolor{purple}{rgb}{0.5,0,0.5} 
\definecolor{turq}{rgb}{0,0.805,0.816} 
\definecolor{maroon}{rgb}{0.51,0,0}
\definecolor{MAROON}{rgb}{0.51,0,0}
\definecolor{redor}{rgb}{0.78,0.078,0.078}
\definecolor{dgreen}{rgb}{0,0.3,0}

\newcommand{\bcdot}{\,{\raise .2ex \hbox{$\centerdot$}}\,}

